\documentclass{article}

\usepackage{amssymb,dsfont,latexsym}

\newtheorem{thm}{Theorem}

\newtheorem{lem}[thm]{Lemma}
\newtheorem{cor}[thm]{Corollary}

\newcommand\enu[1]{\smallskip\newline\makebox[5mm][l]{\rm(#1)}}

\newcommand\bp{\noindent{\it Proof.}\ }

\newcommand\1[1]{{\cal #1}}

\begin{document}

\author{Erling St{\o}rmer}

\date{23-10-2008 }

\title{Duality of cones of positive maps}

\maketitle
\begin{abstract}

We study the so-called K-positive linear maps from $B(L)$ into $B(H)$ for finite dimesional Hilbert 
spaces $L$ and $H$ and give characterizations of the dual cone of the cone of  K-positive maps.  
Applications are given to decomposable maps and their relation to PPT-states.

\end{abstract}

\section*{Introduction}

The study of positive linear maps of C*-algebras, and in particular those of finite dimensions,  has over the last 
years been invigorated by its connection with quantum information theory.  While most work on positive maps
on C*-algebras has been related to completely positive maps, in quantum information theory other classes of maps
appear naturally.  In \cite {st3} the author introduced different cones of positive maps and defined what he
called $K$-positive maps arising from a so-called mapping cone $K$ of positive maps of $B(H)$ into itself, 
see section 2 for details of this and the following. In \cite{BZ}, section 11.2, the authors introduced what they
called the dual cone of a cone of positive maps.  In the present paper we shall follow up this idea by studying
the dual cone of the cone of $K-$positive maps.  Our main result gives several characterizations of when a map
belongs to a dual cone; in particular the result is an extension of the Horodecki Theorem \cite{3 Hor} to
general mapping cones.  Then we show that the dual of a dual cone equals the original cone, and if the mapping
cone $K$ is invariant under the action of the transpose map, then for maps of $B(H)$ into itself the dual cone
consists  of $K^{\sharp}$-positive maps for a mapping cone $K^{\sharp}$ naturally defined by $K.$  Applications are given
to the most studied maps,  like completely positive, copositive, and decomposable maps, and to maps defined
by separable states and PPT-states.  In particular it is shown that if $P$ is the mapping cone of maps which are
both completely positive and copositive, i.e. maps which correspond to PPT-states, then the $P-$positive maps
constitute exactly the dual of the cone of decomposable maps. 

\section{Dual cones}

In this section we shall study certain cones of positive maps from the complex $n\times n$ matrices
$M_{n}$, denoted by $M$ below, to the bounded operators $B(H)$ on a Hilbert space $H,$ 
which we for simplicity assume is finite dimensional.  We denote by $B(M,H)$ (resp. $B(M,H)^+, 
CP(M,H), Cop(M,H))$ the linear (resp. positive, completely positive, and copositive) maps of $M$ into 
$B(H)$. Recall that a map $\phi$ is \textit{copositive} if $t\circ\phi $ is completely
positive, $t$ being the transpose map.  When $M=B(H)$ we shall use the simplified notation
$\1P (H)=B(B(H),H)^+$.  Recall from \cite{st3} that a \textit {mapping cone} is a nonzero
closed cone $K\subset \1 P (H)$ such that if $\phi\in K$ and $a,b\in B(H)$ then the
map $ x\to a\phi(bx b^*)a^*$ belongs to $K.$  Since every completely positive map in
$\1 P(H)$ is of the sum of maps $x\to axa^*$, it follows that if $\phi \in K$ then
$\alpha\circ\phi\circ\beta \in K$ for all $\alpha, \beta \in CP(H)$  - the completely positive
maps in $\1 P(H)$.  We also denote by $Cop(H)$ the copositive maps in $\1 P(H).$ Let
$$
P(M,K) = \{x\in M\otimes B(H): \iota \otimes \alpha(x) \geq 0, \forall \alpha \in K\},
$$
where $\iota$ denotes the identity map. By \cite{st3}, Lemma 2.8, $P(M,K)$ is a proper closed cone.  
If $Tr$ denotes the usual trace on $B(H)$, and also on $M$ and $M\otimes B(H)$ when there is no 
confusion of which algebra we refer to, and $\phi \in B(M,H)^+$, then the dual functional $\tilde \phi$
on $M\otimes B(H)$ is defined by
$$
\tilde \phi (a\otimes b)=Tr(\phi(a) b^t).
$$
We say $\phi$ is \textit{K-positive} if $\tilde \phi$ is positive on the cone $P(M,K).$  It was shown in 
\cite{st3} that if $K=CP(H), $ then $\phi $ is $K-$positive if and only if $\phi$ is completely positive.
Other characterizations will be shown below. We denote by $\1 P_{K}(M)$ the cone of $K-$positive maps.

If $(e_{ij})$ is a complete set of matrix units in $M$ we denote by $p$ the rank 1 operator 
$p=\sum_{ij} e_{ij}\otimes e_{ij} \in M\otimes M,$ and if $\phi\in B(M,H)$ we let $C_{\phi}$ denote
the \textit{Choi matrix} 
$$
C_{\phi}=\iota\otimes\phi (p) = \sum_{ij} e_{ij}\otimes \phi(e_{ij}) \in M\otimes B(H).
$$
Then $\phi$ is completely positive if and only if $C_{\phi}\geq 0$, \cite {C}. We shall show
characterizations of other classes of maps by positivity properties of $C_{\phi}.$  If 
$S\subset B(M,H)^+$ then its \textit{dual cone} is defined by
$$
S^{\flat}= \{\phi\in B(M,H): Tr(C_{\phi}C_{\psi})\geq 0, \forall \psi \in S\}.
$$
If $\phi\in B(M,H)$ we denote by $\phi^t$ the map $t\circ\phi\circ t$. Then $C_{\phi^t}$ is 
the density operator for $\tilde\phi$, i.e. $\tilde\phi(x)=Tr(C_{\phi^t}x)$, see \cite{st5},
Lemma 5.  If $K$ is a mapping cone we put
$$
K^t=\{\phi^t : \phi\in K\}.
$$
Then $K^t$ is also a mapping cone, and is in many cases equal to $K.$ 

We denote by $\phi^*$ the adjoint map of $\phi$ considered as a linear operator of $B(H)$ into
$M$ associated with the Hilbert-Schmidt structure, viz.
$$
Tr(\phi(a)b)=Tr(a\phi^*(b)),  a\in M, b\in B(H).
$$
We can now state the main result of this section.

\begin{thm}\label{thm1}
Let $H$ be a finite dimesional Hilbert space and $K$ a mapping cone in $\1P(H)$.  Let 
$\1P_{K}(M)^\flat$ be the dual cone of the $K-$positve maps $\1P_{K}(M)$.  Let
 $\phi\in B(M,H).$ Then the following conditions are equivalent.
\enu{i} $\phi\in \1P_{K}(M)^\flat$.
\enu{ii} $C_{\phi}\in P(M,K^t)$.
\enu{iii}  $\tilde\phi\circ(\iota\otimes \alpha^*)\geq 0, \forall \alpha \in K$.
\enu{iv} $\alpha\circ\phi \in CP(M,H), \forall \alpha \in K^t$.
\end{thm}
In \cite{st6} we showed a version of the Horodecki Theorem \cite{3 Hor} which in the
case when $n\leq dim H$ is equivalent to the Horodecki Theorem.  We obtain this and 
more as a corollary to Theorem 1, thus showing that it can be viewed as an extension of 
the Horodecki Theorem to arbitrary mapping cones, see also \cite{HSR}. Note that $\1 P(H)$ is a mapping cone 
containing all others, see \cite{st3}, Lemma 2.4.

\begin{cor}\label{cor1}
Let $\phi\in B(M,H)^+$.  Then the following four conditions are equivalent.
\enu{i} $\phi\in \1 P_{\1 P(H)}(M)^{\flat}.$
\enu{ii} $\iota\circ\alpha(C_{\phi})\geq 0, \forall \alpha\in \1P(H).$
\enu{iii} $\tilde\phi \circ (\iota\otimes\alpha) \geq 0,\forall \alpha\in \1 P(H).$
\enu{iv} $\alpha\circ\phi\in CP(M,H), \forall \alpha \in \1 P(H)$.

Furthermore, if $n\leq dim H$ then the above conditions are equivalent to
\enu{v} $ \tilde\phi $ is separable.
\end{cor}
\bp Since $\1 P(H) = \1 P(H)^t$ the four conditions, (i)-(iv) are equivalent by the 
corresponding conditions in Theorem 1.  By \cite{st6}, Lemma 9, (iii)$\Leftrightarrow$ (v)
when $n\leq dim H$, proving the last part of the corollary. 

\medskip
For the proof of the theorem we shall need some lemmas.  The first can easily be extended to
the general situation studied in \cite{st3}.

\begin{lem}\label{lem3}
Let $\rho$ be a linear functional on $M\otimes B(H)$ with density operator $h$.  Let $K$ be a mapping cone in $\1 P(H)$. 
Then $h\in P(M,K)$ if and only if 
$\rho\circ (\iota\otimes \alpha^*) \geq 0, \forall \alpha \in K$.
\end{lem}
\bp
$\rho\circ (\iota\otimes \alpha^*)(x)=Tr(h(\iota\otimes\alpha^*)(x))
=Tr(\iota\otimes\alpha(h)x)$.  Hence $\rho\circ (\iota\otimes \alpha^*)\geq 0$ for all 
$\alpha \in K$ if and only if $\iota\otimes \alpha(h)\geq 0$ for all $\alpha\in K$
if and only if $h\in P(M,K)$, completing the proof.

\medskip
Recall that if $\phi\in B(M,H)^+$ then $\phi^t(x)=\phi(x^t)^t=t\circ\phi\circ  t (x) $.

\begin{lem}\label{lem4}
Let  $\phi\in B(M,H)^+$  Then
\enu{i} $\tilde\phi^{t} ={\tilde\phi }\circ {(t \otimes t)} $.
\enu{ii} $C_{\phi^t}=t\otimes t (C_{\phi})= C_{\phi}^t$.
\enu{iii} $\tilde\phi^t(x)=\tilde\phi(x^t)$.
\end{lem}
\bp
(i) By definitiion of $\phi^t$ and $\tilde\phi$ we have
$$
\tilde\phi^t(a\otimes b)=Tr(\phi^t(a)b^t)=Tr(\phi(a^t)b^{tt})
=\tilde\phi(a^t\otimes b^t)=\tilde\phi\circ (t\otimes t)(a\otimes b).
$$
(ii) By (i) and the fact that $C_{\phi}$ is the density operator for $\tilde\phi^t$, we have
\begin{eqnarray*}
Tr(C_{\phi}a\otimes b)& =&\tilde\phi^t(a\otimes b)=Tr(\phi^t(a)b^t)=Tr(\phi(a^t)^t b^t)\\
&=&
\tilde\phi(a^t \otimes b^t) = Tr(C_{\phi^t} a^t\otimes b^t) = Tr(t\otimes t (C_{\phi^t})a\otimes b).
\end{eqnarray*}
Hence $C_{\phi}= t\otimes t(C_{\phi^t})$, and (ii) follows.

(iii) By (ii) 
$$
\tilde\phi(x^t)=Tr(C_{\phi^t}x^t)=Tr(C_{\phi}^t x^t)=Tr(C_{\phi}x)=\tilde\phi^t (x),
$$
completing the proof.
\begin{lem}\label{lem5}
If $K$ is a mapping cone, then
$$
P(M,K^t)=\{C_{\phi}: C_{\phi}^t \in P(M,K)\}= t\otimes t (P(M,K)).
$$
\end{lem}
\bp  We have
$$
\iota\otimes \alpha^t (x)=(\iota\otimes t)\circ (\iota\otimes \alpha)\circ (\iota\otimes t)(x)=
(t\otimes t)\circ (\iota\otimes\alpha)\circ (t\otimes t)(x).
$$
Thus $x\in P(M,K^t)$ if and only if $t\otimes t (x) \in P(M,K)$, if and only if $x\in t\otimes t (P(M,K))$,
and the two cones are equal.

Each operator $x\in M\otimes B(H)$ is of the form $C_{\phi}$ for some map $\phi \in B(M,H)$.  By Lemma 4(ii)
and the above we therefore have that $C_{\phi^t}=t\otimes t (C_{\phi}) \in P(M,K)$ if and only if 
$C_{\phi}\in t\otimes t (P(M,K)) = P(M,K^t)$, completing the proof.

\medskip
Proof of Theorem 1.

(i)$\Leftrightarrow $(ii) As before let $p=\sum_{ij}e_{ij}\otimes e_{ij},$ where $(e_{ij})$ is a complete 
set of matrix units for $M=M_{n}$.  Then $C_{\phi}=\iota\otimes \phi (p).$  By \cite{st3}, Theorem 3.6,
the cone of $K-$positive maps $\1P_{K}(M)$ is generated as a cone by maps of the form $\alpha\circ \psi$
with $\alpha\in K^d = \{t\circ \alpha^* \circ t: \alpha\in K\}$, and $\psi\in CP(M,H)$. We thus have
 $$
 \phi\in \1P _{K}(M)^{\flat} \Leftrightarrow Tr(C_{\phi}C_{\alpha\circ\psi})\geq 0, \forall \alpha\in K^d,
 \psi\in CP(M,H),
 $$
 hence if and only if
 $$
 0\leq Tr(C_{\phi}( \iota\otimes\alpha )(C_{\psi}))=Tr(\iota\otimes\alpha^* (C_{\phi})C_{\psi}), \forall \alpha, \psi
 $$
 as above, which holds if and only if $\iota\otimes\alpha^* (C_{\phi})\geq 0,  \forall \alpha\in K^d$, since
 $(M\otimes B(H))^+ = \{C_{\psi}: C_{\psi}\in CP(M,H)\}$. Thus $ \phi\in \1P _{K}(M)^{\flat}$ if and
 only if $\iota\otimes\alpha^t (C_{\phi})\geq 0,  \forall\alpha\in K$ if and only if $C_{\phi}\in P(M,K^t)$.
 
 (i)$\Leftrightarrow $(iii) By Lemma 5 $C_{\phi}\in P(M,K^t)$ if and only if $C_{\phi^t}\in P(M,K)$.  Hence by
 (i)$\Leftrightarrow $(ii) $ \phi\in \1P _{K}(M)^{\flat}$ if and only if $C_{\phi^t}\in P(M,K)$, hence by Lemma 3,
 if and only if $\tilde\phi\circ (\iota\otimes\alpha^*)\geq 0, \forall \alpha\in K$, proving (i)$\Leftrightarrow $(iii).
 
 (ii)$\Leftrightarrow $(iv) Recalling the definitions of $C_{\phi}, P(M,K)$ and their properties we have
 $$
 C_{\phi}\in P(M,K^t) \Leftrightarrow  C_{\alpha^t \circ \phi} =\iota\otimes\alpha^t (C_{\phi})\geq 0, \forall 
 \alpha\in K \Leftrightarrow \beta\circ\phi \in CP(M,H), \forall \beta \in K^t.
 $$
 This completes the proof of the theorem.
 
 \medskip
 We conclude the section by showing that taking the dual is a well behaved property.
 \begin{thm}\label{thm5}
 Let $K$ be a mapping cone.  Then $(\1 P_{K}(M)^{\flat})^{\flat}=\1 P_{K}(M)$.
\end{thm}
 \bp
 Let $\phi\in B(M,H)^+$.  Then
 \begin{eqnarray*}
 \phi\in(\1 P_{K}(M)^{\flat})^{\flat}&\Leftrightarrow & Tr(C_{\phi}C_{\psi})\geq 0, \forall \psi\in \1 P_{K}(M)^{\flat} \\
 & \Leftrightarrow &Tr(C_{\phi}C_{\psi})\geq 0, \forall C_{\psi}\in P(M,K^t)\ \  by \ Thm. 1\\
 &\Leftrightarrow &Tr(C_{\phi}( t\otimes t) (C_{\rho})), \forall C_{\rho}\in P(M,K)\ \ by \ Lem. 5\\
 &\Leftrightarrow &Tr(C_{\phi}^t C_{\rho})\geq 0, \forall C_{\rho}\in P(M,K)\\
 &\Leftrightarrow & \tilde\phi(C_{\rho})\geq 0,   \forall C_{\rho}\in P(M,K)\ \ by \ Lem.4(ii)\\
 &\Leftrightarrow & \phi\in \1 P_{K}(M).
 \end{eqnarray*}
 The proof is complete.

 \medskip
 In the notation of \cite{st3} $S(H) $ denotes the mapping cone consisting of maps of the form $x\mapsto 
 \sum \omega_{i}(x)b_{i}$, with $\omega_{i}$ states of $B(H)$ and $b_{i}\in B(H)^+$.  By \cite{st5},
 Theorem 1, $\phi$ is $S(H)-$positive if and only if $\tilde\phi$ is separable.  It follows from Lemma 2.1 in \cite{st3}
 that $\phi$ is positive if and only if $\tilde\phi$ is positive on the cone $B(H)^+\otimes B(H)^+$ generated by
 operators of the form $a\otimes b$ with $a,b\in B(H)^+$, which holds if and only if $Tr(C_{\phi}x)\geq 0$ for
 all $x\in B(H)^+\otimes B(H)^+$.  But by the above $P(B(H),P(H)) = B(H)^+\otimes B(H)^+$.
 Thus $\phi\in \1P(H)$ if and only if $\phi\in \1P_{S(H)}(B(H))^{\flat}$,
 and by Theorem 6 $\phi$ is $S(H)-$positive if and only if $\phi\in \1P(H)^{\flat}$. 

\section{Maps on $B(H)$}

In the previous section we gave characterizations of the dual cone $\1 P_{K}(M)^{\flat}$ of a mapping cone $K$.
A natural question is whether $\1 P_{K}(M)^{\flat}$ equals the cone $\1P_{K^{\sharp}}(M)$ for some mapping cone 
$K^{\sharp}.$  In the present section we shall do this for maps in $\1P(H) =B(B(H),H)^+$ when $K$ is a mapping
cone invariant under the transpose map, viz. $K=K^t.$  The cone $K^{\sharp}$ is defined in our first lemma.

\begin{lem}\label{lem7}
Let $K$ be a mapping cone.  Let $C^{K}$ denote the closed cone generated by all cones
$$
\iota\otimes\alpha^* ((M\otimes B(H))^+),  \alpha\in K.
$$
Let
$$
K^{\sharp}=\{\beta\inÊ\1P(H): \iota\otimes \beta (x) \geq 0,  \forall x\in C^K\}.
$$
Then $K^{\sharp}$ is a mapping cone, and furthermore
$$
K^{\sharp}=\{\beta\inÊ\1P(H): \beta\circ\alpha^* \in CP(H), \forall \alpha\in K\}.
$$
\end{lem}
\bp
Let $\beta\in K^{\sharp}$ and $\gamma\in CP(H)$. Then clearly $\gamma\circ\beta\in K^{\sharp}$.
If $\alpha \in K$ and $\gamma\in CP(H)$, then
$$
(\beta\circ\gamma)\circ\alpha^* = \beta\circ(\gamma\circ\alpha^*) = \beta\circ(\alpha\circ\gamma^*)^*.
$$
Since $K$ is a mapping cone, and $\gamma\in CP(H), \alpha\circ\gamma^*\in K$, hence
$\iota\otimes\beta\circ(\alpha\circ\gamma^*)^*\geq 0,$ so that $\beta\circ \gamma\in K^{\sharp},$ proving
the first part of the lemma.  To show the second part we have that $\beta\in K^{\sharp}$ if and only if
$$
\iota\otimes \beta\circ\alpha^*(x) = \iota\otimes \beta (\iota\otimes\alpha^*)(x)\geq 0, \forall 
x\geq 0, \alpha\in K,
$$
which holds if and only if $\beta\circ\alpha^* \in CP(H)$, because a map $\gamma\in \1P(H)$ is completely positive
if and only if $\iota\otimes \gamma \geq 0$ on $B(H\otimes H)^+$. The proof is complete. 

\medskip
We shall need the following rephrasing of Choi «s result  \cite{C} that $\phi$ is completely positive 
if and only if $C_{\phi}\geq 0.$

\begin{lem}\label{lem8}
Let $\phi\in \1P(H)$ and $\omega$ the maximally entangled state, viz. $\omega(x)= \frac{1}{n}Tr(px),
n=dim H. $ Then $\phi\in CP(H)$ if and only if $\omega\circ (\iota\otimes\phi)\geq 0.$
\end {lem}
\bp
$\phi$ is completely positive if and only if
$$
0\leq Tr(C_{\phi}x)= Tr(p (\iota\otimes\phi^*)(x)) = n \omega\circ (\iota\otimes\phi^*)(x), \forall x\geq 0.
$$
Since $CP(H)$ is closed under the *-operation, the lemma follows.

\begin{lem}\label{lem9}
Let $K$ be a mapping cone and $C^K$ as in Lemma 7. Let $\phi\in \1P(H)$. Then we have
\enu{i} $\tilde\phi$ is positive on $C^K$ if and only if $\phi^{*t}\in K^{\sharp}$.
\enu{ii} $\tilde\phi^t$ is positive on $C^K$ if and only if $\phi^*\in K^{\sharp}$.

If $K=K^t$ then  $\tilde\phi$ is positive on $C^K$ if and only if $\tilde\phi^t$ is positive on $C^K$. 
\end{lem}
\bp
(i) Let $\alpha\in K.$  Then
\begin{eqnarray*}
\tilde\phi(\iota\otimes\alpha^*(x))&=&Tr(C_{\phi^t}(\iota\otimes\alpha^*)(x))\\
&=& Tr(\iota\otimes\phi^t(p)(\iota\otimes\alpha^*)(x))\\
&=& Tr(p( \iota\otimes (\phi^{t*}\circ\alpha^*))(x))\\
&=& Tr(p (\iota\otimes( \phi^{*t}\circ\alpha^*))(x)),
\end{eqnarray*}
since $\phi^{t*}=\phi^{*t}.$  By Lemma 8 it follows that $\tilde\phi{\geq 0}$ on $C^K$ if and only if 
 $ \phi^{*t}\circ\alpha^*\in CP(H)$, hence if and only if $\phi^{*t}\in K^{\sharp}$, by Lemma 7.
 
 (ii) follows from (i) by applying (i) to $\phi^t.$
 
 Assume $K=K^t$.  If $\phi^{*t}\in K^{\sharp}$ then by Lemma 7 $\phi^{*t}\circ\alpha^* \in CP(H)$ 
 for all $\alpha\in K$, hence by assumption for all $\alpha^{*t}, \alpha\in K.$ Thus
 $(\phi^*\circ\alpha^*)^t = \phi^{*t}\circ\alpha^{*t} \in CP(H),$ hence $\phi^*\circ\alpha^*
 \in CP(H)$ for all $\alpha\in K.$ Thus $\phi^*\in K^{\sharp}$, so by (i) and (ii), if $\tilde\phi\geq 0$ on
 $C^K$ then $\tilde\phi^t\geq 0$ on $C^K.$  Similarly we get the converse implication. The proof is complete.
 
 \begin{lem}\label{lem10}
 Let $\pi\colon B(H)\otimes B(H)\to B(H)$ be defined by $\pi(a\otimes b)=b^t a.$  Then the function
 $Tr\circ \pi$ is positive and linear. Furthermore, if $\phi\in \1P(H)$ then
 $$
 \tilde\phi = Tr\circ\pi\circ (\iota\otimes\phi^{*t}).
 $$
 \end{lem}
 \bp
 Linearity is clear. To show positivity let $x=\sum a_{i}\otimes b_{i}\in B(H)\otimes B(H). $ Then
 \begin{eqnarray*}
 Tr(xx^*)&=&  \sum Tr\circ\pi (a_{i}a_{j}^*\otimes b_{i}b_{j}^*)\\
 &=& \sum Tr(b_{j}^{t*}b_{i}^t  a_{i} a_{j}^*)\\
 &=& \sum Tr(a_{j}^* b_{j}^{t*} b_{i}^t a_{i})\\
 &=& Tr((\sum b_{j}^t a_{j})^*(\sum b_{i}^t a_{i}))\geq 0,
 \end{eqnarray*}
 so $Tr\circ\pi$ is positive.  The last formula follows from the computation
 $$
 \tilde\phi(a\otimes b) = Tr(\phi(a)b^t) = Tr(a\phi^*(b^t))
 = Tr(a\phi^{*t}(b)^t) = Tr\circ\pi (\iota\otimes \phi^{*t}(a\otimes b)).
 $$
 
 \begin{lem}\label{lem11}
 Assume $K=K^t$ for $K$ a mapping cone. Then $C^K = P(B(H),K^{\sharp})$.
 \end{lem}
 \bp
 By definition of $K^{\sharp}, C^K\subset P(B(H),K^{\sharp}).$  Suppose
 $y_{0}\in B(H)\otimes B(H)$ and $y_{0}$ is not in $C^K$.  By the Hahn-Banach 
 Theorem there is a linear functional $\tilde\phi$ on $B(H)\otimes B(H)$ which is
 positive on $C^K$, and $\tilde\phi(y_{0})< 0.$  By Lemma 9 and the assumption 
 that $K=K^t, \tilde\phi^t \geq 0$ on $C^K$ as well, and $\phi^{*t}\in K^{\sharp}.$
 Write $y_{0}$ in the form $y_{0}=\sum a_{i}\otimes b_{i}, a_{i},b_{i}\in B(H),$ 
 and let $\pi$ be as in Lemma 10.  then we have
 \begin{eqnarray*}
 Tr\circ \pi(\iota\otimes \phi^{*t}(y_{0}))&=& Tr\circ\pi(\sum a_{i}\otimes \phi^*(b_{i}^t)^t)\\
 &=& Tr(\sum \phi^*(b_{i}^t)a_{i})\\
 &=& Tr(\sum b_{i}^t \phi(a_{i})) = \tilde\phi(y_{0}) < 0.
 \end{eqnarray*}
 Since by Lemma 10, $Tr\circ\pi$ is positive, $\iota\otimes \phi^{*t}(y_{0})$ is not a
 positive operator, hence $y_{0} $ does not belong to $P(B(H),K^{\sharp}),$ hence
 $C^K = P(B(H),K^{\sharp})$.  The proof is complete.
 
 \medskip
 We can now prove the main result in this section, which shows that every map in $\1 P_{K}(B(H))^{\flat}$
 is $K^{\sharp}-$positive when $K=K^t$.  Thus Theorem 1 yields equivalent conditions for $K-$positivity
 when $K=K^t$, by replacing $K^{\sharp}$ by $K$,
 \begin{thm}\label{thm12}
 Let $K$ be a mapping cone such that $K=K^t$.  Then $\1P_{K}(B(H))^{\flat} =\1P_{K^{\sharp}}(B(H))$,
 so that  $\1P_{K}(B(H)) =\1P_{K^{\sharp}}(B(H))^{\flat}$.
 \end{thm}
 \bp
 By Theorem 1(iii) $\phi\in \1P_{K}(B(H))^{\flat}$ if and only if $\tilde\phi$ is positive on $C^K,$ 
 hence by Lemma 11 if and only if $\phi$ is $K^{\sharp}-$positive, i.e. $\phi\in \1P_{K^{\sharp}}(B(H)).$
 The last statement folows from Theorem 6.
 
 \section{Decomposable maps}
 
 It was shown in \cite{3 Hor} that a state $\rho$ on $M_{2}\otimes M_{2}$ or $M_{2}\otimes M_{3}$ is 
 separable if and only if $\rho\circ (\iota\otimes t)\geq 0$, i.e. if and only if it is a PPT state (equivalently, $\rho$
 is said to satisfy the Peres condition).  They did this by using the fact that $B(M_{2},\0C^2)^+$ and $B(M_{2},\0C^3)^+$
 consist of decomposable maps, i.e. maps which are sums of completely positive and copositive maps. We shall in
 this section show characterizations of decomposable maps which yield characterizations of  PPT states.
 
 Let  $D$ denote the set of $\alpha\in \1P(H)$ such that $ \alpha$ is decomposable.  Let $P$ denote the set of
  $\alpha\in \1P(H)$ such that $ \alpha$ is both completely positive and copositive. Thus
 $$
 D= CP(H) \vee Cop(H). \ \ P=CP(H) \cap Cop(H).
 $$
 
 \begin{thm}\label{thm13}
 Let $M=M_{n}$, and  $D$ and $P$ as above. Then $D$ and $P$ are mapping cones and satisfy the identities
 $$
 \1P_{P}(M)^{\flat}=\1P_{D}(M),\ \  \1P_{P}(M)= \1P_{D}(M)^{\flat}.
 $$
\end{thm} 

 \medskip
 Note that by \cite{st5},Proposition 4, $P$ consists of the maps $\phi\in \1P(H)$ such that $\tilde\phi$
 is a PPT state. The proof of Theorem 13 is divided into some lemmas.
 
 \begin{lem}\label{lem14}
 In the notation of Theorem 13 let $E$ and $F$ denote the cones
 $$
 E=\{x\in M\otimes B(H): x\geq 0, \ \ or  \ \iota\otimes t(x)\geq 0\}.
 $$
 $$
 F=\{x\in M\otimes B(H): x\geq 0, \ \ and  \ \iota\otimes t(x)\geq 0\}.
 $$
 Then the following conditions are equivalent for $\phi\in B(M,H)^+.$
 \enu{i} $\phi$ is both completely positive and copositive.
 \enu{ii} $C_{\phi}\in F$.
 \enu{iii} $\tilde\phi \geq 0$ on $E.$
 \end{lem}
 \bp
 Note that $E$ and $F$ are closed under the action of $\iota\otimes t$.
 
 (i)$\Rightarrow$(ii) Since $\phi$ is completely positive $C_{\phi}\geq 0$, and since $\phi$ is
 copositive $\iota\otimes t(C_{\phi})\geq 0.$  Thus $C_{\phi}\in F$.
 
 (ii)$\Rightarrow$ (i) If $x\geq 0$ then $\tilde\phi(x)=Tr(C_{\phi^t}x)\geq 0$, so $\phi$
 is completely positive. If $\iota\otimes t(x)\geq 0$ then
 $$
 \tilde\phi(x)=Tr(C_{\phi^t}x) = Tr(\iota\otimes t (C_{\phi^t})\iota\otimes t (x))\geq 0,
 $$
 so $\phi$ is copositive.
 
 (ii)$\Rightarrow$ (iii) The same argument as for  (ii)$\Rightarrow$ (i) applies.
 
 (iii)$\Rightarrow$ (ii) If $x\geq 0$ then $Tr(C_{\phi^t}x)= \tilde\phi(x)\geq 0$.  Similarly, if 
  $\iota\otimes t(x)\geq 0$ then
  $$
  0\leq \tilde\phi(x)= Tr(C_{\phi^t}x)=Tr(\iota\otimes t (C_{\phi^t})\iota\otimes t (x))\geq 0,
  $$
  hence $\iota\otimes t (C_{\phi^t})\geq 0.$ The proof is complete.
  
  \begin{lem}\label{lem15}
  Let $K$ and $L$ be mapping cones in $\1P(H)$ and $M$ as before. Then
  $$
  P(M,K\vee L)=P(M,K)\cap P(M,L).
  $$
  \end{lem}
  \bp
  $x\in P(M,K\vee L)$ if and only if $\iota\otimes \alpha(x)\geq 0$ for all $ \alpha\in K\cup L$,
  if and only if  $\iota\otimes \alpha(x)\geq 0$ whenever $\alpha\in K$ or $\alpha\in L$, if and only if 
  $x\in P(M,K)\cap P(M,L),$ proving the lemma.
  
  \medskip
  Note that we did not use that $M$ is the $n\times n$ matrices in the above proof, so the lemma is true
   for $M$ replaced by an operator system.
   
   \medskip
  
     \begin{lem}\label{lem17}
     Let $E$ be as in Lemma 14. Then
     $$
     P(M,P)=P(M,CP(H))\vee P(M,Cop(H)) =E.
     $$
     \end{lem}
     \bp
     Since $P(M,CP(H))=(M\otimes B( H))^+$ and $P(M,Cop(H))=\iota\otimes t (M\otimes B(H))^+)$, it
     is clear that $P(B(M,CP(H))\vee P(M,Cop(H)) =E.$
     
     If $K$ and $L$ are mapping cones with $K\subset L$ and $M$ is a finite dimensional $C^*-$algebra then clearly 
     $P(M,L)\subset P(M,K)$, so clearly $P(B(H),P)$ contains the right side of the lemma.  Suppose the inclusion is proper.
     By the Hahn-Banach Theorem there exists a linear functional $\tilde\phi$ which is positive on $E$ and for some
     $x\in P(B(H),P), \tilde\phi(x)< 0.$  By Lemma 14 $\phi$  is both completely positive and copositive, hence so
     is $\phi^{*t},$ so that by Lemma 10 
     $$
     \tilde\phi(x)= Tr\circ\pi (\iota\otimes \phi^{*t}(x))\geq 0,
     $$
     a contradiction. This proves  equality of the cones.
     
     \begin{lem}\label{lem19}
     With the previous notation we have
     $$
      P(M,D)=F=\{C_{\beta}: \beta\in\1P_{P}(M)\}.
     $$
     $$
     \1P_{P}(M)=\{\beta:C_{\beta}\in F\}.
     $$ 
     \end{lem}
     \bp
     By \cite{st3}, Theorem 3.6, $\1P_{P}(M) $ is generated by maps $\alpha\circ\psi$ with $\alpha\in P, \psi\in CP(M,H).$ 
     By Lemma 15 
     $$
     F=P(M,CP(H))\cap P(M,Cop(H)) = P(M,D).
     $$
     Let $\gamma\in D.$ Then, with $\alpha$ and $\psi$ as above, 
     $$
     \iota\otimes\gamma (C_{\alpha\circ\psi}) = \iota\otimes\gamma\circ \iota\otimes\alpha(C_{\psi})=
    \iota\otimes\gamma\circ\alpha(C_{\psi})\geq 0,
    $$
    because   $\psi\in CP(M,H),$ so  $C_{\psi}\geq 0,$ and $\gamma\circ\alpha$ is completely positive since $\gamma$ is
    decomposable and $\alpha\in P.$
    Thus $C_{\alpha\circ\psi}\in P(M,D) = F.$  It follows that $ C_{\beta}\in P(M,D)=F$ for all $\beta\in \1P_{P}(M).$
    Since $D=D^t$, by Lemma 5
    $$
    P(M,D)=\{C_{\beta}:C_{\beta}\in P(M,D)\} =\{C_{\beta}:C_{\beta}\in F\}.
    $$
    By Lemma 16 if $\beta\in B(M,H)$ then $\beta\in \1P_{P}(M)$ if and only if $\tilde\beta \geq 0$ on $E,$ so by Lemma 14,
    if and only if $ C_{\beta}\in F.$  Hence $\1P_{P}(M) =\{\beta: C_{\beta}\in F\}, $ hence
    $$
    P(M,D)=F =\{C_{\beta}:\beta\in \1P_{P}(M)\},
    $$
    as asserted. The proof is complete.
    
    \medskip
    Proof of Theorem 13.
    
    This is immediate, since by Lemma 18 $ \phi\in \1P_{P}(M)^{\flat}$ if and only if $\tilde\phi(C_{\beta})\geq 0$
    for all $\beta\in \1P_{P}(M)$ if and only if $\tilde\phi\geq 0$ on $P(M,D),$ i.e. $\phi\in \1P_{D}(M).$  The 
    other identity follows from the first and Theorem 6.
    
    \medskip
    For the rest of this section we consider the case when $M=B(H).$

     \begin{thm}\label{thm19}
     Let $D$ and $P$ be as in Theorem 13. Then $D=P^{\sharp}.$ Furthermore $\1P_{P}(B(H))^{\flat}=\1P_{D}(B(H))=D.$
     \end{thm}
     \bp
     To simplify notation let $\1P_{P}=\1P_{P}(B(H))$  and similarly for $D.$ By Theorem 13 $\1P_{P}^{\flat}=
     \1P_{D},$ and by Theorem 12  $\1P_{P}^{\flat}=\1P_{P^{\sharp}}.$ Thus $\1P_{D}=\1P_{P^{\sharp}}.$
     Note that by Lemma 7
     $$
     P^{\sharp}=\{\beta\in \1P(H): \beta\circ\alpha^* \in CP(H), \forall \alpha\in P\}.
     $$
     hence $D\subset P^{\sharp}$.  Since  $\1P_{D}=\1P_{P^{\sharp}}$, a linear functional $\tilde\phi$ is
     positive on $P(B(H),D)$ if and only if it is positive on $P(B(H),P^{\sharp})$.  Since  $D\subset P^{\sharp}$
     it follows that  $P(B(H),P^{\sharp})\subset P(B(H),D),$ hence from the Hahn-Banach Theorem that they are
      equal.  Thus by Lemma 9 and Lemma 11 $\phi^*\in P^{\sharp}$ if and only if $\tilde\phi\geq 0 $ on
      $C^P=P(B(H),P^{\sharp})=P(B(H),D)$, if and only if $\phi\in \1P_{D} =\1P_{P}^{\flat}$.
      
      Let  $\phi^*\in P^{\sharp}$, so that $\phi\in \1P_{D} =\1P_{P}^{\flat}$, hence by Theorem 1 $C_{\phi}
      \in P(B(H),P)$, which by Lemma 15 equals the cone $E$ in Lemma 14.  Thus $C_{\phi}=C_{\phi_{1}}+
      C_{\phi_{2}}$ with $\phi_{1}\in CP(H)$ and $C_{\phi_{2}}\in Cop(H),$ hence $\phi\in D,$ as is $\phi^*.$
     Thus $ P^{\sharp}\subset D$, and they are equal.  Since $D$ is closed under the *-operation, so is   $ P^{\sharp},$
     so by the above $D=P^{\sharp}=\1P_{D}$. The proof is complete.
     
     \medskip
     We conclude by showing the analogue of Lemma 8 for decomposable maps.  When $M=B(H)  $ the result is a 
     strenghtening of the result in \cite{st2}, which states that $\phi$ is decomposable if and only if $\iota\otimes\phi$ is
     positive on $F.$
     \begin{cor}\label{cor20}
     Let $\omega$ denote the maximally entangled state on $B(H)$, and let $\phi\in \1P(H).$ Then $\phi$ is
     decomposable if and only if $\omega\circ(\iota\otimes \phi)\geq 0$ on the cone $F=\{x\in B(H\otimes H)^+
     :\iota\otimes t (x)\geq 0\}$.
     \end{cor}
     \bp
     By Theorem 18 and Lemma 16 $\phi\in D$ if and only if $\phi\in \1P_{D}$ if and only if $\tilde\phi\geq 0$ on
     $P(M,D)=F,$ if and only if 
     $$
     0 \leq Tr(C_{\phi}x)=Tr(\iota\otimes\phi(p)x)=Tr(p (\iota\otimes\phi^*)(x)) = \omega\circ(\iota\otimes\phi^*)(x),
     $$
     for all $x\in F.$ Since $D$ is closed under *-operation, the corollary follows.

Department of Mathematics, University of Oslo, 0316 Oslo, Norway.

e-mail erlings@math.uio.no


\begin{thebibliography}{999}



 \bibitem{BZ}
 I.Bengtsson and K. Zyczkowski, {\em Geometry of quantum states}, Cambridge Universty Press (2006).
\bibitem{C}
M-D. Choi, {\em Completely positive linear maps on complex matrices}, Linear alg. and its applic. 10,
(1975), 285-290.
 \bibitem{3 Hor}
 M.Horodecki, P.Horodecki, and R.Horodecki, {\em Separability of
 mixed states: necessary and sufficient conditions}, Physics
 Letters, A 223, (1996), 1-8.
\bibitem{HSR}
M.Horodecki, P.W.Shor, and M.B.Ruskai, {\em Entanglement breaking channels}, Rev. Mat. Phys. 15,
(2003), 629-641.
\bibitem{st2}
E.St{\o}rmer, {\em Decomposable positive maps on $C^*$-algebras}, Proc.
Amer. Math. Soc. 86 (1982), 402-404.
 
 \bibitem{st3}
 E.St{\o}rmer, {\em Extension of positive maps into $B(H)$}, J.
 Funct. Anal. 66, No.2 (1986), 235-254.


 \bibitem{st5}
 E.St{\o}rmer, {\em Separable states and positive maps},
 J.Funct.Anal. 254 (2008), 2303-2312.
\bibitem{st6}
E.St{\o}rmer, {\em Separable states and positive maps II}, Math.
Scand. to appear.
 
 \end{thebibliography}
\end{document}